\documentclass [10pt,twoside]{article}
\usepackage[T1]{fontenc}
\usepackage{multicol}
\usepackage{amsfonts,amsmath,amssymb,amsthm,bm}
\usepackage[english]{babel}
\usepackage{fancybox}
\usepackage{graphicx}
\usepackage[utf8]{inputenc}
\usepackage{enumerate}
\usepackage{fullpage}
\usepackage{dsfont}
\usepackage{stmaryrd}
\usepackage{color}
\usepackage{float}

\usepackage{tikz}
\usepackage{pgfplots}

\numberwithin{equation}{section}

\date\today

\usepackage{geometry}
\usepackage[french]{layout}
\usepackage{fancyhdr} 
 \usepackage{amsopn}
\DeclareMathOperator{\Supp}{Supp}
\DeclareMathOperator{\Div}{div}

\newcommand{\mb}[1]{\mathbb{#1}}
\newcommand{\mc}[1]{\mathcal{#1}}

\theoremstyle{definition} 
\newtheorem{cond}{\textsc{Condition}}[section]

\theoremstyle{remark}
\newtheorem{rem}{\textsf{Remark}}

\theoremstyle{plain}
\newtheorem{theo}{\textsc{Theorem}}
\newtheorem{Lemme}{\textsf{Lemma}}[section]
\newtheorem{prop}{\textsc{Proposition}}[section]
\newtheorem{ass}{\textsf{Assumption}}[section]

\newcommand{\cqfd}
{%
\mbox{}%
\nolinebreak%
\hfill%
\rule{2mm}{2mm}%
\medbreak%
\par%
}

\author{{\normalsize Michel Duprez}\footnote{
Institut de Math\'ematiques de Marseille (I2M),
Centre de Math\'ematiques et Informatique (CMI),
Technop\^ole Ch\^ateau-Gombert, Bureau 211, 39, rue F. Joliot Curie
13453 Marseille Cedex 13, France\texttt{michel.duprez@univ-fcomte.fr}, 
},~{\normalsize Pierre Lissy}\footnote{CEREMADE, Universit\'e Paris-Dauphine \& CNRS UMR 7534, PSL, 75016 Paris, France, \texttt{lissy@ceremade.dauphine.fr.}}}

\title{Positive and negative results on the internal controllability of parabolic equations coupled by zero and first order terms}

\date{\today}

\begin{document}

\maketitle

\begin{abstract}
 This paper is devoted to studying the null and approximate controllability 
  of two linear  coupled parabolic equations posed on a smooth domain $\Omega$ of $\mathbb R^N$ ($N\geqslant 1$)  
  with coupling terms of zero and first orders and one control
 localized in some arbitrary nonempty open subset  $\omega$ of the domain $\Omega$. 
We prove the null controllability under a new sufficient condition 
and we also provide 
the first example of a not approximately controllable system 
in the case where the support of one of the nontrivial first order coupling terms intersects the control domain $\omega$.

\end{abstract}


\vspace{0.3cm}



\section{Introduction}
\subsection{Presentation of the problem and main results}

\hspace*{4mm} Let $T>0$,  let $\Omega$ be a bounded domain of $\mb{R}^N$ ($N\in\mb{N}^*$) 
 of class $\mathcal C^2$
and  let $\omega$ be an arbitrary nonempty open subset of $\Omega$. Let $Q_T:=(0,T)\times\Omega$, 
$q_T:=(0,T)\times\omega$ and $\Sigma_T:=(0,T)\times\partial\Omega$. 
We consider the following system of two parabolic linear equations with variable coefficients and coupling terms of order zero and one 
\begin{equation}\label{system primmal}
 \left\{\begin{array}{ll}
\partial_ty_1=\Div (d_1\nabla y_1)+g_{11}\cdot\nabla y_1+g_{12}\cdot\nabla y_2+a_{11}y_1+a_{12}y_2+\mathds{1}_{\omega}u&\mbox{in }Q_T,\\\noalign{\smallskip}
\partial_ty_2=\Div (d_2\nabla y_2)+g_{21}\cdot\nabla y_1+g_{22}\cdot\nabla y_2+a_{21}y_1+a_{22}y_2&\mbox{in } Q_T,\\\noalign{\smallskip}
y=0&\mbox{on }\Sigma_T,\\\noalign{\smallskip}
y(0,\cdot)=y^0&\mbox{in }\Omega,
        \end{array}
\right.
\end{equation}
where $y^0\in L^2(\Omega)^2$ is the initial condition and $u\in L^2(Q_T)$ is the control.

The  zero and first  order coupling terms
  $(a_{ij})_{1\leqslant i,j\leqslant 2}$ and $(g_{ij})_{1\leqslant i,j\leqslant 2}$ 
  are assumed (for the moment) to be  in  $ L^\infty(Q_T)$ and 
 in $L^\infty(0,T;W^{1}_{\infty}(\Omega))^N$, respectively.
For $l\in\{1,2\}$, the second order elliptic self-adjoint operator $\Div (d_l\nabla)$ is given by
 \begin{equation*}
  \Div (d_l\nabla)=\sum\limits_{i,j=1}^N\partial_i(d^{ij}_l\partial_j),
 \end{equation*}
with 
\begin{equation*}
\left\{
 \begin{array}{l}
  d^{ij}_l\in W^{1}_{\infty}(Q_T),\\\noalign{\smallskip}
   d^{ij}_l= d^{ji}_l\mbox{ in }Q_T,
 \end{array}\right.
\end{equation*}
where the coefficients $d^{ij}_l$ satisfy the uniform ellipticity condition 
\begin{equation*}
\sum\limits_{i,j=1}^N d^{ij}_l\xi_i\xi_j\geqslant d_0|\xi|^2\mbox{ in }Q_T,~\forall \xi\in\mb{R}^N,
\end{equation*}
for a constant $d_0>0$.
 
 It is well-known (see for instance \cite[Th. 3, p. 356-358]{MR2597943}) 
that for  every initial data  $y^0\in L^2(\Omega)^2$ and every control $u\in L^2(Q_T)$, 
 System \eqref{system primmal} 
 admits a unique solution $y$  
in 
$W(0,T)^2$, where  
\begin{equation*}
W(0,T):=L^2(0,T;H^1_0(\Omega ))
\cap H^1(0,T;H^{-1}(\Omega ))\hookrightarrow\mc{C}^0([0,T];L^2(\Omega)).
\end{equation*}

In this article, we are concerned with the approximate or null controllability of System \eqref{system primmal}. Let us recall the precise definitions of these notions. We say that System  \eqref{system primmal} is 
\begin{itemize}
\item[$\bullet$]\textit{approximately controllable} on the time interval $(0,T)$ if for every initial condition $y^0\in L^2(\Omega)^2$, every target  $y^T\in L^2(\Omega)^2$
and every $\varepsilon>0$, there exists a control $u\in L^2(Q_T)$ such that the corresponding solution $y$ to System \eqref{system primmal} satisfies
\begin{equation*}
\|y(T,\cdot)-y^T\|_{L^2(\Omega)^2}\leqslant \varepsilon.
\end{equation*} 
\item[$\bullet$]\textit{null controllable} on the time interval $(0,T)$ if for every initial condition $y^0\in L^2(\Omega)^2$, 
there exists a control $u\in L^2(Q_T)$ such that the corresponding solution $y$ to System \eqref{system primmal} satisfies
\begin{equation*}
y(T,\cdot)=0\mbox{ in }\Omega.
\end{equation*} 
\end{itemize}

It is well-known that if a parabolic system like \eqref{system primmal}  is null controllable on the time interval $(0,T)$, 
 then it is also approximately controllable on the time interval $(0,T)$
 (this is an easy consequence of usual results 
 of backward uniqueness for parabolic equations as given for example in \cite{MR0338517}).

Since the case $a_{21}\neq0$ and $g_{21}=0$ in $(t_0,t_1)\times\omega_0\subset q_T$ 
has already been studied in \cite{gonzalez2010controllability},  we will always work under the following assumption.
\begin{ass}\label{MainAss}
There  exists  $t_0<t_1\in (0,T)$ and  a nonempty open subset $\omega_0$ of $\omega$  such that 
\begin{equation*}
g_{21}\neq0\mbox{ in }(t_0,t_1)\times\omega_0. 
\end{equation*}
\end{ass}

Moreover, as we will see in Section \ref{sec:simpl coupl}, it is possible, with the help of  appropriate changes of variables and unknowns (we lose a little bit of regularity on the coefficients though, see Section \ref{sec:simpl coupl}), 
to replace the coupling operator 
$g_{21}\cdot \nabla +a_{21}$
 by the simpler coupling operator $\partial_{x_1}$ (where $x_1$ is the first direction in space), at least locally on some subset of $q_T$.

Hence, without loss of generality, we can also work under the following assumption.

\begin{ass}\label{MainAss2}
There  exists a nonempty open subset $ \mathcal O_T$ of 
$\omega_0$  such that 
$$g_{21}\cdot \nabla +a_{21}=\partial_{x_1}\mbox{ on } \mathcal O_T:=(t_0,t_1)\times \mathcal O .$$
\end{ass}

For a nonempty set $\omega_T\subset \mathbb{R}^{N+1}$, 
let us denote by $\mc{C}^0_{t,x_2,...,x_N}( \overline{\omega}_T)$ the subset of $\mc{C}^0(\overline{\omega}_T)$ 
composed by the functions depending only on the variables $t,~x_2,~x_3,...,x_N$. Let us now introduce the following condition, which will be crucial in our following results, and  which is closely related to the particular form for the coupling term given in Assumption \ref{MainAss2} (removing this assumption would make Condition \ref{cond:modul} impossible to write down explicitly).

\begin{cond}\label{cond:modul}
There exists a nonempty open set $\omega_T\subset (t_0,t_1)\times\mathcal O$ such that
\begin{equation}\label{alg resol:cond cont}
\left\{\begin{array}{l}
\widetilde{a}_{22}\mbox{ is not an element of the }
\mc{C}^0_{t,x_2,...,x_N}(\overline{\omega}_T)\mbox{-module }\\\noalign{\smallskip}
\left\langle 1,\widetilde{g}_{22}^2,...,\widetilde{g}_{22}^N,d_2^{22},...,d_2^{NN}\right\rangle_{\mc{C}^0_{t,x_2,...,x_N}(\overline{\omega}_T)},
\end{array}\right.
\end{equation}
where 
\begin{equation}\label{resol alg:def g22 tilde}
\left\{\begin{array}{l}
\widetilde{g}_{22}^i:=g_{22}^i-\sum\limits_{j=1}^N\partial_{x_j}d_{22}^{ij},\\\noalign{\smallskip}
\widetilde{a}_{22}:=- a_{22}+\Div (g_{22}).
\end{array}\right.
\end{equation}
\end{cond}

Our first main result  is the following:

\begin{theo}\label{theo:positif}
Let $d_i^{kl}, g_{ij}^k\in \mc{C}^{N^2+3}(\overline{\omega}_T)$ and $a_{ij}\in \mc{C}^{N^2+2}(\overline{\omega}_T)$ 
for every $i,j\in\{1,2\}$ and $k,l\in\{1,...,N\}$.
Assume that Assumptions \ref{MainAss}, \ref{MainAss2}  and 
Condition \ref{cond:modul} hold.
Then System \eqref{system primmal} is null controllable at any time $T>0$.
\end{theo}

\begin{rem}
Theorem \ref{theo:positif} is stated and will be proved in the case of two coupled parabolic equations and one control. However, as in \cite{ML15}, it is possible to extend Theorem \ref{theo:positif} 
to  systems of $m$ parabolic  equations controlled by $m-1$ controls for arbitrary $m\geqslant 2$.
More precisely, consider the system
\begin{equation}\label{system primmalg}
 \left\{\begin{array}{ll}
\partial_ty_1=\Div(d_1\nabla y_1)+\sum_{i=1}^mg_{1i}\cdot\nabla y_i+\sum_{i=1}^{m}a_{1i}y_i+\mathds{1}_{\omega}u_1&\mbox{in }Q_T,\\\noalign{\smallskip}
\partial_ty_2=\Div(d_2\nabla y_2)+\sum_{i=1}^mg_{2i}\cdot\nabla y_i+\sum_{i=1}^{m}a_{2i}y_i+\mathds{1}_{\omega}u_2&\mbox{in }Q_T,\\\noalign{\smallskip}
\vdots\\\noalign{\smallskip}
\partial_ty_{m-1}=\Div(d_{m-1}\nabla y_{m-1})+\sum_{i=1}^mg_{(m-1)i}\cdot\nabla y_i+\sum_{i=1}^{m}a_{(m-1)i}y_i+\mathds{1}_{\omega}u_{m-1}&\mbox{in } Q_T,\\\noalign{\smallskip}
\partial_ty_m=\Div (d_m\nabla y_m)+\sum_{i=1}^mg_{mi}\cdot\nabla y_i+\sum_{i=1}^{m}a_{mi}y_i&\mbox{in }Q_T,\\\noalign{\smallskip}
y_1=\ldots=y_m=0&\mbox{on } \Sigma_T,\\\noalign{\smallskip}
y_1(0,\cdot)=y_1^0,\ldots,~y_m(0,\cdot)=y_m^0&\mbox{in } \Omega,
        \end{array}
\right.
\end{equation}
where $y^0:=(y_1^0,\ldots,y_m^0)\in L^2(\Omega)^m$ is the initial data 
and $u:=(u_1,\ldots,u_{m-1})\in L^2(Q_T)^{m-1}$ is the control. 
Let us suppose that there exists $i\in \{1,...,m\}$,  $t_0<t_1\in (0,T)$ 
and  a nonempty open subset $\omega_0$ of $\omega$  such that $g_{mi}(t,x)\neq0$ on $q_T:=(t_0,t_1)\times\omega_0$. 
As explained in Section \ref{sec:simpl coupl}, we can suppose that the operator 
$g_{mi}\cdot\nabla +a_{mi}$ is equal to $\partial_{x_1}$ in $(t_0,t_1)\times\mathcal O\subset q_T$. 
Assume that there exists an open set $\omega_T\subset (t_0,t_1)\times\mathcal O$ such that
\begin{equation*}
\left\{\begin{array}{l}
\widetilde{a}_{mm}\mbox{ is not an element of the }
\mc{C}^0_{t,x_2,...,x_N}(\overline{\omega}_T)\mbox{-module }\\\noalign{\smallskip}
\left\langle 1,\widetilde{g}_{mm}^2,...,\widetilde{g}_{mm}^N,d_2^{22},...,d_2^{NN}\right\rangle
_{\mc{C}^0_{t,x_2,...,x_N}(\overline{\omega}_T)},
\end{array}\right.
\end{equation*}
where 
\begin{equation*}
\left\{\begin{array}{l}
\widetilde{g}_{mm}^i:=g_{mm}^i-\sum\limits_{j=1}^N\partial_{x_j}d_{mm}^{ij},\\\noalign{\smallskip}
\widetilde{a}_{mm}:=- a_{mm}+\Div (g_{mm}).
\end{array}\right.
\end{equation*}
Then we can adapt the proof of Theorem \ref{theo:positif} to prove that  System \eqref{system primmalg} is null controllable on the time interval $(0,T)$ under suitable regularity conditions on the coefficients.
\end{rem}
\begin{rem}Condition \ref{cond:modul} is clearly technical since it does not even cover the case of constant coefficients proved in \cite{ML15}, the general case given in \cite{benabdallah2014} (under some assumption on the control domain) or the one-dimensional result given in \cite{duprez2016}. 
However, Theorem \ref{theo:negatif} implies that one  cannot expect the null controllability to be true in general without extra assumptions on the coefficients. We do not know what would be a reasonable necessary and sufficient condition on the coupling terms for the null controllability of System \eqref{system primmal}.
\end{rem}

The second main result of the present paper is the following surprising result.

\begin{theo}\label{theo:negatif2}
Let us assume that $\omega\subset\subset \Omega$.
Let $\omega_1$  be a nonempty regular open set satisfying $\omega\subset\subset \omega_1\subset\subset\Omega$.
and consider a function $\theta\in \mathcal{C}^{\infty}(\overline{\Omega})$ satisfying
\begin{equation*}
\left\{\begin{array}{l}
\theta=1\mbox{ in }\omega,\\\noalign{\smallskip}
\Supp(\theta)\subset\overline\omega,\\\noalign{\smallskip}
\theta>0\mbox{ in } \omega_1
\end{array}\right.
\end{equation*}
Then there exists $a\in \mc{C}^{\infty}(\overline{\Omega})$ such that the system 
\begin{equation}\label{syst simpl 1D}
\left\{\begin{array}{ll}
\partial_ty_1=\Delta y+\mathds{1}_{\omega}u&\mbox{in } Q_T,\\\noalign{\smallskip}
\partial_ty_2=\Delta y_2+ay_2+\partial_{x_1}(\theta y_1)&\mbox{in } Q_T,\\\noalign{\smallskip}
y=0&\mbox{on } \Sigma_T,\\\noalign{\smallskip}
y(0,\cdot)=y^0&\mbox{in }\Omega
\end{array}\right.
\end{equation}
is not approximately controllable (hence not null controllable) on the time interval $(0,T)$.
\end{theo}

In other words, Theorem \ref{theo:negatif2} tells us that for every  control set $\omega$ strongly included in $\Omega$, there exists a  potential $a$ for which approximate controllability of \eqref{syst simpl 1D} does not hold, in any space dimension. We may improve a bit this result on the one-dimensional case, where we are able to obtain the following result, which expresses that for some well-constructed potential $a$,  that there exists one control domain on which System \eqref{syst simpl} is not approximately controllable (hence not null controllable) and another control domain on which System \eqref{syst simpl} is null controllable (hence approximatively controllable),
 highlighting the surprising fact that some \emph{geometrical conditions} on the control domain has to be imposed in order to obtain a controllability result.

\begin{theo}\label{theo:negatif}
Consider the following system
\begin{equation}\label{syst simpl}
\left\{\begin{array}{ll}
\partial_ty_1=\partial_{xx} y+\mathds{1}_{\omega}u&\mbox{in } (0,T)\times(0,\pi),\\\noalign{\smallskip}
\partial_ty_2=\partial_{xx} y_2+ay_2+\partial_{x} y_1&\mbox{in } (0,T)\times(0,\pi),\\\noalign{\smallskip}
y(\cdot,0)=y(\cdot,\pi)=0&\mbox{on } (0,T),\\\noalign{\smallskip}
y(0,\cdot)=y^0&\mbox{in }(0,\pi).
\end{array}\right.
\end{equation}

There exists a coefficient $a\in \mc{C}^{\infty}([0,\pi])$ such that:
\begin{enumerate}
\item There exists an open interval $(a,b)\subset\subset (0,\pi)$
such that, for all $T>0$, System \eqref{syst simpl} is null controllable (then approximatively controllable) 
at time $T$.
\item There exists an open interval  $(a,b)\subset\subset (0,\pi)$
such that, for all $T>0$, System \eqref{syst simpl} is not approximatively controllable (then not null controllable) 
at  time $T$.
\end{enumerate}
\end{theo}

\begin{rem}
Let us mention that Theorems \ref{theo:negatif2} and \ref{theo:negatif} are the first negative result for the controllability of  System \eqref{system primmal} 
when the support of the first order coupling term intersects the control domain in the case of distributed controls.
 The authors want to highlight that the coupling operator is constant in the whole domain
and  nevertheless  the system can be controllable or not following the localisation of the control domain,  which is an unexpected phenomenon.
\end{rem}


 \subsection{State of the art}

\hspace*{4mm}  
Many
models of interest involve  (linear or non-linear) coupled equations of parabolic systems, notably in fluid mechanics, medicine, chemistry, ecology, geology, etc., and this explains why
during the past years, the study of the controllability properties of linear or nonlinear parabolic systems has been an increasing subject  of interest (see for example the survey \cite{ammar2011recent}). The main issue is what is  called the \emph{indirect} controllability, that is to say one wants to control many equations with less controls than equations, by acting indirectly on 
the equations where no control term appears thanks to the coupling terms appearing in the system. This notion is fondamental for real-life applications, since in some complex systems only some quantities can be effectively controlled.
Here, we will concentrate on the previous results concerning the null or approximate controllability  of linear parabolic systems with distributed controls, but 
 there are also many other results concerning  boundary controls or other classes of 
 systems like hyperbolic systems.

  First of all, in the case of zero order coupling terms, the case of constant coefficients is now completely treated and we refer to  \cite{ammar2009generalization} and \cite{ammar2009kalman}   for parabolic systems having constant coupling coefficients (with diffusion coefficients that may depend on the space variable though) and for some results in the case of time-dependent coefficients. In the case of zero and one order coupling terms and constant coefficients, a necessary 
and sufficient condition in the case of $m$ equations and $m-1$ controls for constant coefficients is provided in  \cite{ML15} by the authors.

The case of space-varying coefficients remains still widely open despite many new partial results these last years. In the case where the support of the coupling terms intersects
the control domain, a general result is proved in \cite{gonzalez2010controllability} 
for parabolic systems in cascade form with one control force (and possibly one order coupling terms). 
We also mention  \cite{MR2226005}, where a result of null controllability is proved in 
the case of a system of two equations with one control force, with an application to the controllability 
of a nonlinear system of transport-diffusion equations.
In the situation where the coupling regions do not intersect the control domain, 
the situation is still not very well-understood and we have partial results, in general under technical and geometrical restrictions, notably on the control domain
(see for example \cite{alabau2013}, \cite{MR3039207}, \cite{MR2783322} and  \cite{CherifMinimalTimeDisjoint}). Let us mention that in this case, there might appear a minimal time for the null controllability of System \eqref{system primmal}
(see \cite{Ammar-Khodja2015}), which is a very surprising phenomenon for parabolic equations, because of the infinite speed of propagation of the information.

Concerning the case of first order coupling terms, we mention \cite{gonzalez2010controllability} which gives some controllability results 
when the coefficient $g_{21}$ is equal to zero on the control domain.
Let us also mention  the recent work \cite{benabdallah2014}, which concerns the small systems in small dimension, that is to say   $2\times2$ and $3\times 3$ systems. 
The authors of \cite{benabdallah2014} suppose that the control domain contains a part of the boundary $\partial\Omega$.
Recently, in \cite{duprez2016}, the first author studied a particular cascade system with space dependent coefficients and in dimension one thanks to the moment method, and obtained necessary and sufficient conditions on the coupling terms of order $0$ and $1$ for the null controllability.
To conclude, let us also mention another result given in  \cite{ML15} by the authors, which provides a sufficient condition for null controllability in dimension one  for space and time-varying coefficients  under some  technical conditions on the coefficients, which turns out to be exactly equivalent to Condition \ref{cond:modul} under Assumption \ref{MainAss2} (but with more regularity than in Assumption \ref{MainAss}). Hence,  Theorem \ref{theo:positif} can be seen as a generalization in the multi-dimensional case of the one-dimensional result given in \cite{ML15}.
For a more detailed state of the art concerning this problem, 
we refer to \cite{ML15}.

 Hence, the present paper improves the previous results in the following sense:
\begin{itemize}
\item Contrary to \cite{benabdallah2014,guerrerosyst22,duprez2016,ML15}, 
we prove in Theorem \ref{theo:positif} the null controllability of System  \eqref{system primmal} 
with a condition on $a_{22}$ but for space/time dependent coefficients, in any space dimension 
and without any condition on the control  domain.

\item In the previous results, it was surprising to have some very different sufficient conditions for the null controllability 
of System  \eqref{system primmal} in the case of first order coupling terms, for example on one hand 
constant coupling coefficients and on the other hand a region of control which intersects the boundary of the domain. 
Through the example of a not approximately controllable system given in Theorem \ref{theo:negatif2} and \ref{theo:negatif}, 
we can now better understand why such different conditions appeared since the expected general condition for the null controllability of System \eqref{system primmal} with space and time-varying coefficients (i.e. it is sufficient that the control and coupling region intersect) may be false in general if $\omega\subset\subset\Omega$.
\end{itemize}

\section{Simplification of the coupling term}\label{sec:simpl coupl}
In this section, we will prove that it is possible to replace locally the coupling operator $g_{21}\cdot \nabla +a_{21}$ by $\partial_{x_1}$, 
where $x_1$ is the first direction in space. This kind of simplification has already been used in \cite[Lemma 2.6]{benabdallah2014} for example, and we refer to this article for a more detailed proof (see also \cite{duprez2016}). Let us remark that the regularities stated in Lemma \ref{reduc1} are higher than the one stated in Theorem \ref{theo:positif} due to technical reasons appearing in the proofs of Lemmas \ref{reduc1} and \ref{reduc2}.

\begin{Lemme}\label{reduc1}
Let $ d_i^{kl},~g_{ij}^k,~a_{ij}\in \mc{C}^{N^2+4}([t_0,t_1]\times\overline{\omega}_0)$
for every $i,j\in\{1,2\}$ and $k,l\in\{1,...,N\}$.
Suppose that Assumption \ref{MainAss} is verified. Then,
there exist a nonempty open subset $U$ of $\mathbb R^{N-1}$, a positive real 
number $\varepsilon$ and a  
$\mathcal{C}^{N^2+3}$-diffeomorphism $\Lambda$ from $U_{\varepsilon}:=(t_0,t_1)\times(0,\varepsilon)\times U $  
to an open set 
$(t_0,t_1)\times \mathcal O\subset (t_0,t_1)\times\omega_0$ that keeps $t$ invariant and such that if we call 
$\widetilde y_1:=y_1\circ\Lambda$ and $\widetilde y_2:=y_2\circ\Lambda$, 
then there exist a matrix $\widetilde d_2\in \mathcal M_N(\mathcal C^{N^2+3}(U_{\varepsilon}))$, a vector $\widetilde g_{22}\in (\mathcal C^{N^2+3}(U_{\varepsilon}))^N$ and coefficients $\widetilde a_{21},~\widetilde a_{22}\in \mathcal C^{N^2+3}(U_{\varepsilon})$ such that locally on $U_{\varepsilon}$ one has 

\begin{gather}\label{imp}
\partial_t\widetilde y_2=\Div (\widetilde d_2\nabla \widetilde y_2)+\widetilde g_{22}\cdot\nabla  \widetilde y_2+\widetilde a_{22} \widetilde y_2 
+\partial_{x_1} \widetilde{y_1}+\widetilde a_{21} \widetilde y_1\mbox{ in } U_{\varepsilon}.
\end{gather}
\end{Lemme}


\textbf{Proof of Lemma~\ref{reduc1}}\\
Let us consider some open hyper-surface $\gamma$ of class $\mathcal{C}^\infty$ included in $\omega_0$ on which $g_{21}\cdot\nu<0$, 
where $\nu$ is the normalized outward normal on $\gamma$ 
(this can always be done since $g_{21}\not =0$ on $(t_0,t_1)\times\omega_0$ and is at least continuous), 
small enough such that it can be parametrized by a local  diffeomorphism
$$F:s_0:=(s_2,\ldots,s_N)\in U\subset\mathbb R^{N-1}\mapsto F(s_0)\in\gamma,$$
where $U$ is a nonempty open set.
We call $\gamma_T:=(t_0,t_1)\times\gamma$.
Let us consider some $\mathcal{C}^{N^2+4}$ extension of $g_{21}$ (that exists thanks to the regularity of $\gamma$ and $g_{21}$) that we denote by $g^T_{21}:(t,x)\in\mathbb R^{N+1}\mapsto(0,g_{21}(t,x))\in\mathbb R^{N+1}$. 
Using the Cauchy-Lipschitz Theorem, we infer that for every $(t,\sigma)\in \gamma_T$, there exists a unique global solution to the Cauchy Problem
 \begin{equation*}
 \left\{\begin{array}{ll}
\frac{d}{ds}\Phi(t,s,\sigma)=g^T_{21}(\Phi(t,s,\sigma)),\\\noalign{\smallskip}
\Phi(t,0,\sigma)=(t,\sigma).
        \end{array}
\right.
\end{equation*}
Since $\Phi$ is continuous and $g_{21}\cdot\nu<0$ on $\gamma_T$,
we deduce that there exists some $\varepsilon>0$ such that $\Phi(t,s,\sigma)\in (t_0,t_1)\times\omega_0$ for every $s\in (0,\varepsilon)$ and every $(t,\sigma)\in \gamma_T$.
We define 
$$\Lambda:(t,s,z)\in (t_0,t_1)\times(0,\varepsilon)\times U\mapsto \Phi(t,s,F(z)).$$
Then, by the inverse mapping theorem, $\Lambda$ is a $\mathcal C^{N^2+4}$-diffeomorphism from $U_{\varepsilon}$ 
to  $(t_0,t_1)\times\mathcal O:=\Lambda(U_{\varepsilon})$ with $\mathcal O\subset\omega_0$. 
Let us call $\widetilde y_1(t,s,z):=y_1(\Lambda(t,s,z))$ and $\widetilde y_2(t,s,z):=y_2(\Lambda(t,s,z))$, then it is clear that 
$$\partial_t \widetilde y_i(t,s,z)=(\partial_t y_i)\circ\Lambda(t,s,z) \mbox{ for }i=1,2~ \mbox{ and }\partial_s \widetilde y_2(t,s,z)=(g_{21}\cdot\nabla y_2)\circ\Lambda(t,s,z),$$
and hence we obtain \eqref{imp} and the regularities wished for the new coefficients by writing down the equation verified by $\widetilde y$.
\cqfd

Let us now perform a second useful reduction.
\begin{Lemme}\label{reduc2}
There exists an open subset$\mc{O}_T$ of $ U_{\varepsilon}$ and 
a function $\theta\in \mathcal C^{N^2+4}(\Omega)$ such that $|\theta(x)|\geqslant C$ for some constant $C>0$ and if 
$$\overline y_1(t,x):=\theta^{-1}(t,x) \widetilde y_1(t,x)$$
and $$\overline y_2(t,x):=\theta^{-1}(t,x) \widetilde y_2(t,x),$$
then  there exists some coefficients $\overline a_{22}\in \mathcal C^{N^2+2}(\mc{O}_T)$ and $\overline g_{22}\in \mathcal C^{N^2+3}(\mc{O}_T)^N$ such that locally on $\mathcal O_T$ one has

\begin{gather}\label{imp2}
\partial_t\overline  y_2=\Div (\widetilde d_2\nabla \overline y_2)+\partial_{x_1}\overline y_1+\overline g_{22}\cdot\nabla \overline y_2
+\overline a_{22} \overline y_2 \mbox{ in } {\mathcal O_T}.
\end{gather}
\end{Lemme}

\textbf{Proof of Lemma~\ref{reduc2}}\\
Let us consider some $\theta\in \mathcal C^{N^2+4}(\overline{\Omega})$ such that $|\theta(x)|\geqslant C$ for some constant $C>0$, and consider the change of unknowns
$$\left\{\begin{array}{l}
\overline y_1(t,x):=\theta^{-1}(t,x) \widetilde y_1(t,x),\\\noalign{\smallskip}
\overline y_2(t,x):=\theta^{-1}(t,x) \widetilde y_2(t,x).\end{array}\right.$$
Using equation \eqref{imp}, we infer that $\overline y_2$ verifies 
$$\partial_t\overline  y_2=\Div (\widetilde{d}_2\nabla \overline y_2)+\overline g_{22}\cdot\nabla \overline y_2+\overline a_{22} \overline y_2
+\partial_{x_1} \overline y_1 +\theta^{-1}(\partial_{x_1}\theta+\widetilde{a}_{21}\theta)\overline y_1,$$
where $\overline g_{22}:=2\theta^{-1}\widetilde{d}_2\nabla\theta+\widetilde{g}_{22}$ and $\overline a_{22}:=\theta^{-1}\Div(\widetilde{d}_2\nabla\theta)+\theta^{-1}\widetilde{g}_{22}\nabla\theta+\widetilde{a}_{22}$. 
Hence, if we choose $\theta\in \mathcal C^{N^2+4}(\overline{\Omega})$ satisfying $\partial_{x_1}\theta+\widetilde{a}_{21}\theta=0$ and $|\theta(x)|\geqslant C$ 
in $Q_T$, which is always possible, 
then  $\overline y_1$ and $\overline y_2$  verify \eqref{imp2} and we have $\overline a_{22}\in \mathcal C^{N^2+2}(\mc{O}_T)$ and $\overline g_{22}\in \mathcal C^{N^2+3}(\mc{O}_T)^N$ .

\cqfd

\section{Proof of Theorem \ref{theo:positif}} 
During all this Section, we will always assume that Assumptions \ref{MainAss} and \ref{MainAss2} are satisfied.
\subsection{Strategy : Fictitious control method }\label{section strategy}

The fictitious control method has already been used for instance  in \cite{gonzalezperez2006}, \cite{coronlissy2014}, \cite{ACO}, \cite{CG16} and \cite{ML15}. 
Roughly, the method is the following: we first control the equations with two controls (one on each equation) 
and we try to eliminate the control on the last equation thanks to algebraic manipulations locally on the control domain. 
For more details, see for example \cite[Section 1.3]{ML15}.
Let us be more precise and decompose the problem into three different steps:

\vspace*{0.2cm}

\begin{itemize}

\item[(i)]
\textbf{Analytic Problem: Null controllability by two forces}\\
Find a solution  $(\widehat{y},\widehat{u})$ in an appropriate space 
to the control problem by two controls 
\begin{equation}\label{strat:syst lin pb ana}
 \left\{\begin{array}{ll}
\partial_t\widehat{y}_1=\Div (d_1\nabla \widehat{y}_1)+g_{11}\cdot\nabla \widehat{y}_1+g_{12}\cdot\nabla \widehat{y}_2+a_{11}\widehat{y}_1+a_{12}\widehat{y}_2+
\widehat{u}_1&\mbox{in }Q_T,\\\noalign{\smallskip}
\partial_t\widehat{y}_2=\Div (d_2\nabla \widehat{y}_2)+g_{21}\cdot\nabla \widehat{y}_1+g_{22}\cdot\nabla \widehat{y}_2+a_{21}\widehat{y}_1+a_{22}\widehat{y}_2+\widehat{u}_2&\mbox{in } Q_T,\\\noalign{\smallskip}
\widehat{y}=0&\mbox{on }\Sigma_T,\\\noalign{\smallskip}
\widehat{y}(0,\cdot)=y^0,~\widehat{y}(T,\cdot)=0&\mbox{in }\Omega,
        \end{array}
\right.
\end{equation}
where the controls $\widehat{u}_1$ and $\widehat{u}_2$ are regular enough 
and with a  support  strongly included in $\omega_T$ (remind that $\omega_T$ was introduced in Condition \ref{cond:modul}).
Solving Problem \eqref{strat:syst lin pb ana} is easier than solving the null controllability 
on the time interval $(0,T)$ of System \eqref{system primmal}, 
because we control System \eqref{strat:syst lin pb ana}  with one control on each equation.  
The important point is that the control 
has to be regular enough, so that it can be differentiated a certain amount 
of times with respect to the space and/or time variables (see the next section about the algebraic resolution). 
\begin{prop}\label{prop:contr regul}
Let $k\in\mathbb{N}^*$. Suppose that 
$ d_i^{kl},~g_{ij}^k\in \mc{C}^{k+2}(\overline{\omega}_T)$ and $a_{ij}\in \mc{C}^{k+1}(\overline{\omega}_T)$
for every $i,j\in\{1,2\}$ and $k,l\in\{1,...,N\}$.
Then there exists two constants $K>0$ and $C_k$ 
such that for every initial condition $y^0\in L^2(\Omega)^2$ 
one can find a control $u\in \mathcal{C}^k(Q_T)^2$  verifying moreover 
$\Supp(u)\subset \subset \omega_T$
for which 
the solution to System \eqref{strat:syst lin pb ana} is equal to zero at time $T$ 
and the following estimate holds:
\begin{equation}\label{intro:exp control}
\|u\|_{\mathcal{C}^k(Q_T)^2}\leqslant C_k\|y^0\|_{L^2(\Omega)^2}.
\end{equation}
\end{prop}

The controllability of parabolic systems with regular controls is nowadays well-known. 
For a proof of Proposition \ref{prop:contr regul}, one can adapt the
strategy developed in  
\cite{gonz_perez_insensit04,perez_gonz_insent04,bodart_burgos_perez_local_04,
gonzalezperez2006}
where the authors prove the controllability of parabolic systems with $L^{\infty}$ controls
thanks to the fictitious control method and the local regularity of parabolic equations.
For more details, we refer to \cite[Chap. I, Sec. 2.4]{these_michel}.
It is also possible to use the Carleman estimates 
(see for instance \cite{MR1751309} and \cite[Section 2.3]{ML15}), 
however this will impose the coefficients of System \eqref{strat:syst lin pb ana}  to be regular in the whole space $Q_T$ (and would require higher regularity on $\Omega$).



\item[(ii)]
\textbf{Algebraic Problem: Null controllability  by one force}\\
For given $\widehat{u}_1,\widehat{u}_2$ with $\Supp(\widehat{u}_1,\widehat{u}_2)$ strictly included 
in $\omega_T$, find $(z,v)$,  in an appropriate space, satisfying the following control problem:
\begin{equation}\label{strat:probleme ramene a tout l espace}
 \left\{\begin{array}{ll}
\partial_tz_1=\Div (d_1\nabla z_1)+g_{11}\cdot\nabla z_1+g_{12}\cdot\nabla z_2+a_{11}z_1+a_{12}z_2+\widehat{u}_1+v&\mbox{in }\omega_T,\\\noalign{\smallskip}
\partial_tz_2=\Div (d_2\nabla z_2)+\partial_{x_1} z_1+g_{22}\cdot\nabla z_2+a_{22}z_2+\widehat{u}_2&\mbox{in } \omega_T,
        \end{array}
\right.
\end{equation}
with $\Supp(z,v)$ strictly included in $\omega_T$, which impose the initial and final  data and the boundary conditions. 
We recall that $g_{21}\cdot\nabla+a_{21}$ is equal to $\partial_{x_1}$ in $\omega_T$. 
 We will solve this problem using the notion of \emph{algebraic resolvability} 
of differential systems, which is based on ideas coming from \cite[Section 2.3.8]{Gromovbook} 
and was already used in some different contexts in \cite{coronlissy2014}, 
\cite{ACO}, \cite{ML15} or \cite{CG16}.
The idea is to write System 
\eqref{strat:probleme ramene a tout l espace} as an \emph{underdetermined} 
 system in the variables $z$ and $v$ and to see $\widehat{u}$ as a source term. 
More precisely, we remark that System \eqref{strat:probleme ramene a tout l espace} can be rewritten as
 \begin{equation}\label{def L Strategy 2}
 \mathcal{L}(z,v)=f,
 \end{equation}
 where $f:=\widehat{u}$ and
 \begin{equation*}
\mathcal{L}(z,v):= 
 \left(\begin{array}{c}
\partial_tz_1-\Div (d_1\nabla z_1)-g_{11}\cdot\nabla z_1-g_{12}\cdot\nabla z_2-a_{11}z_1-a_{12}z_2-v\\\noalign{\smallskip}
\partial_tz_2-\Div (d_2\nabla z_2)-\partial_{x_1} z_1-g_{22}\cdot\nabla z_2-a_{22}z_2
        \end{array}
\right).
\end{equation*}
The goal in Section \ref{sec:resol alg} will be then to find a partial differential operator $\mc{M}$ satisfying
\begin{gather}
\label{strat:LMB}
\mc{L}\circ\mc{M}=Id\mbox{ in }\omega_T.
\end{gather}
Thus to solve control problem \eqref{strat:probleme ramene a tout l espace}, 
it suffices to take 
\begin{equation*}
(z,v):=\mathcal{M}(f).
\end{equation*}
When \eqref{strat:LMB} is satisfied, we say that  System 
\eqref{def L Strategy 2}
is \emph{algebraically solvable}.

\item[(iii)]
\textbf{Conclusion}\\
If we are able to solve the analytic and algebraic problems, 
then it is easy to check that   $(y,u):=(\widehat{y}-z,-v)$ will be a solution to 
System (\ref{system primmal}) 
in an appropriate space and will satisfy $y(T,\cdot)\equiv0$ in $\Omega$ 
(for more explanations, see \cite[Prop. 1]{coronlissy2014} and the proof of Theorem \ref{theo:positif} on pages 11-12).

\end{itemize}

\subsection{Algebraic solvability of the linear control problem}\label{sec:resol alg}

The goal of this section is to solve algebraic problem \eqref{def L Strategy 2}. 
We will use the following lemma:

\begin{Lemme}\label{lemme:resol alg}
Let $\omega$ be a nonempty open subset of $\mathbb{R}^n$ \emph{(}$n\geqslant 1$\emph{)} and let $R\in\mathbb N^*$. 
Consider two differential operators $\mathcal L_1$ and $\mathcal L_2$ defined for every $\varphi\in \mathcal C^{\infty}(\overline{\omega})$ by
\begin{equation*}
 \mathcal{L}_1\varphi:=\partial_{x_1}\varphi
 \mbox{ and }
 \mathcal{L}_2\varphi:=a_0\varphi+\sum\limits_{i=1}^R a_iD^{\alpha_i}\varphi,
\end{equation*}
where, for $\alpha_i=(\alpha_i^2,...,\alpha_i^n)$, $D^{\alpha_i}:=\partial_{x_2}^{\alpha_i^2}\cdots\partial_{x_n}^{\alpha_i^n}$.
If $a_i\in\mc{C}^M(\overline{\omega})$ for every $i\in\{0,...,R\}$ 
where $$M:=\sum\limits_{j=1}^R\beta_j\mbox{ with }\beta_j\mbox{ the order of the operator }\sum\limits_{i=j}^Ra_iD^{\alpha_i}$$
and $a_0$ is not an element of the $\mc{C}^0_{x_2,...,x_n}(\overline{\widetilde\omega})$-module
\begin{equation}\label{module}
\begin{array}{l}
\left\langle a_1,...,a_R\right\rangle_{\mc{C}^0_{x_2,...,x_n}(\overline{\widetilde\omega})},
\end{array}
\end{equation}
for a nonempty open subset $\widetilde\omega$ of $\omega$,
then there exists two differential operators $\mathcal M_1$ and $\mathcal M_2$ such that
\begin{equation}\label{M1M2}
\mathcal M_1\circ \mathcal L_1+\mathcal M_2\circ \mathcal L_2=Id\mbox{ in }\mathcal C^{\infty}(\overline{\widetilde\omega}).
\end{equation}
\end{Lemme}

\textbf{Proof of Lemma \ref{lemme:resol alg}}\\
The goal is to apply some differential operators  $\mathcal M_1$ and $\mathcal M_2$ to $\mathcal{L}_1\varphi$ and $ \mathcal{L}_2\varphi$ in order to obtain $\varphi$.
So, since $\varphi$ is not appearing in $\mathcal{L}_1\varphi$, we would like to eliminate all the derivatives $D^{\alpha_i}\varphi$ in the expression of $\mathcal L_2\varphi$ by differentiations and linear combinations.

If $a_0\not =0$ and $a_i=0$ in $\omega$  for every $i\in \{1,....,R\}$, we define $$\mathcal N:=\mathcal L_2.$$ 
If not, let $k_1$ be the smallest number of $\{1,....,R\}$ such that there exists a nonempty open subset $\omega_1$ of $\omega$ where $|a_{k_1}|>\delta>0$. 
Then we consider $\mc{L}_3$ the commutator of $\mc{L}_1$ and $a_{k_1}^{-1}\mc{L}_2$:
\begin{equation*}
\mc{L}_3\varphi:=[\mc{L}_1,a_{k_1}^{-1}\mc{L}_2]\varphi
=\partial_{x_1}\left(\frac{a_0}{a_{k_1}}\right)\varphi+\sum\limits_{i={k_1}+1}^R \partial_{x_1}\left(\frac{a_i}{a_{k_1}}\right)D^{\alpha_i}\varphi.
\end{equation*}
Again, if for every $i\in \{k_1+1,....,R\}$, we have $\partial_{x_1}\left(\frac{a_i}{a_{k_1}}\right)=0$ in $\omega$, 
we define $$\mathcal N:=\mathcal L_3.$$
If not, let $k_2$ be the smallest number of $\{k_1+1,....,R\}$ such that there exists a nonempty open subset $\omega_2$ of $\omega_1$ 
where $|\partial_{x_1}\left(\frac{a_{k_2}}{a_{k_1}}\right)|>\delta>0$. 
Then we consider $\mc{L}_4$ the commutator of $\mc{L}_1$ and $\left[\partial_{x_1}\left(\frac{a_{k_2}}{a_{k_1}}\right)\right]^{-1}\mc{L}_3$:
\begin{equation*}
\mc{L}_4\varphi:=[\mc{L}_1,\left[\partial_{x_1}\left(\frac{a_{k_2}}{a_{k_1}}\right)\right]^{-1}\mc{L}_3]\varphi
=\partial_{x_1}\left(\frac{\partial_{x_1}\left(\frac{a_0}{a_{k_1}}\right)}{\partial_{x_1}\left(\frac{a_{k_2}}{a_{k_1}}\right)}\right)\varphi
+\sum\limits_{i=k_2+1}^R \partial_{x_1}\left(\frac{\partial_{x_1}\left(\frac{a_i}{a_{k_1}}\right)}{\partial_{x_1}\left(\frac{a_{k_2}}{a_{k_1}}\right)}\right)D^{\alpha_i}\varphi.
\end{equation*}

Again, if, for every $i\in \{k_2+1,....,R\}$, we have $ \partial_{x_1}\left(\frac{\partial_{x_1}\left(\frac{a_i}{a_{k_1}}\right)}{\partial_{x_1}\left(\frac{a_{k_2}}{a_{k_1}}\right)}\right)=0$ in $\omega_2$, 
we define $$\mathcal N:=\mathcal L_4.$$
If not, we continue the same reasoning that will stop at some point since there is only a finite order of derivatives $R$. Hence, we obtain, for a $m\in\{1,...,R\}$, a nonempty open subset $\widetilde\omega$ of $\omega$ and an operator
\begin{equation}\label{expression N}
\mathcal N\varphi
:=\mathcal L_{m+2}\varphi
=\partial_{x_1}\left(\frac{\partial_{x_1}\left(\cdots\frac{\partial_{x_1}\left(\frac{a_{0}}{a_{k_1}}\right)}
{\vdots}\right)}
{\partial_{x_1}\left(\cdots\frac{\partial_{x_1}\left(\frac{a_{k_m}}{a_{k_1}}\right)}
{\vdots}\right)}\right)
\varphi\mbox{ in }\widetilde\omega.
\end{equation}
Moreover, $\mathcal N$ is obtained by making iterated commutators of operators involving only $\mathcal L_1$ and $\mathcal L_2$. Hence it is clear that there exists two linear partial differential operators $\widetilde {\mathcal M}_1$ and $\widetilde {\mathcal M}_2$ such that 
$$\mathcal N=\widetilde {\mathcal M}_1\mathcal L_1+\widetilde {\mathcal M}_2\mathcal L_2.$$
Hence, in view of \eqref{expression N}, we will have the desired conclusion as soon as 
the coefficient in the right-hand side in \eqref{expression N}  is different from zero. Let us explain into more details what this condition exactly means.
For the sake of clarity, let us assume that $m=3$ (but the following reasoning can be extended to any $m\in \{1,\ldots, R\}$). We remark that
\begin{equation}\label{express0}
\partial_{x_1}\left(\frac{\partial_{x_1}\left(\frac{\partial_{x_1}\left(\frac{a_{0}}{a_{k_1}}\right)}
{\partial_{x_1}\left(\frac{a_{k_2}}{a_{k_1}}\right)}\right)}
{\partial_{x_1}\left(\frac{\partial_{x_1}\left(\frac{a_{k_3}}{a_{k_1}}\right)}
{\partial_{x_1}\left(\frac{a_{k_2}}{a_{k_1}}\right)}\right)}\right)
= 0
\end{equation}
holds only if, for some $\lambda_3\in \mathcal C^0_{x_2,...,x_n}(\overline{\widetilde\omega})$, we have
\begin{equation*}
\frac{\partial_{x_1}\left(\frac{\partial_{x_1}\left(\frac{a_{0}}{a_{k_1}}\right)}
{\partial_{x_1}\left(\frac{a_{k_2}}{a_{k_1}}\right)}\right)}
{\partial_{x_1}\left(\frac{\partial_{x_1}\left(\frac{a_{k_3}}{a_{k_1}}\right)}
{\partial_{x_1}\left(\frac{a_{k_2}}{a_{k_1}}\right)}\right)}
= \lambda_3.
\end{equation*}
The last expression can be rewritten as
\begin{equation}\label{expres1}
\partial_{x_1}\left(\frac{\partial_{x_1}\left(\frac{a_{0}-\lambda_3a_{k_3}}{a_{k_1}}\right)}
{\partial_{x_1}\left(\frac{a_{k_2}}{a_{k_1}}\right)}\right)
=0.
\end{equation}
Again, \eqref{expres1} holds only if, for some $\lambda_2\in \mathcal C^0_{x_2,...,x_n}(\overline{\widetilde\omega})$, we have
\begin{equation*}
\frac{\partial_{x_1}\left(\frac{a_{0}-\lambda_3a_{k_3}}{a_{k_1}}\right)}
{\partial_{x_1}\left(\frac{a_{k_2}}{a_{k_1}}\right)}
=\lambda_2,
\end{equation*}
or, equivalently,
\begin{equation*}
\partial_{x_1}\left(\frac{a_{0}-\lambda_3a_{k_3}-\lambda_2a_{k_2}}{a_{k_1}}\right)
=0.
\end{equation*}
Thus \eqref{express0} is satisfied only if, for some $\lambda_1,~\lambda_2,~\lambda_3\in \mathcal C^0_{x_2,...,x_n}(\overline{\widetilde\omega})$, we have
\begin{equation*}
a_{0}=\lambda_3a_{k_3}+\lambda_2a_{k_2}+\lambda_1a_{k_1}.
\end{equation*}
Hence, we find back condition \eqref{cond:modul} and the proof of Lemma \ref{lemme:resol alg} is achieved.

\cqfd

We are now able to prove the algebraic solvability of \eqref{def L Strategy 2}.
\begin{prop}\label{prop:resol alg}
Suppose that $d_i^{kl},~ g_{ij}^k,~a_{ij}\in \mc{C}^{N^2}(\overline{\omega}_T)$ 
for every $i,j\in\{1,2\}$ and $k,l\in\{1,...,N\}$. 
Then, under Condition \ref{cond:modul}, System \eqref{def L Strategy 2}  is algebraically solvable 
with an operator $\mc M$  of order $N^2$.
\end{prop}

\textbf{Proof of Proposition \ref{prop:resol alg}}\\
Let us remark that the first equation of System  \eqref{def L Strategy 2} can be rewritten locally on $\omega_T$ as
$$v=\partial_tz_1-\Div (d_1\nabla  z_1)-g_{11}\cdot\nabla  z_1-g_{12}\cdot\nabla  z_2-a_{11}  z_1-a_{12}  z_2-f_1,$$
hence one can always solve algebraically first the second equation of System \eqref{def L Strategy 2}, 
 $v$ will then be given with respect to $z_1$, $z_2$ and $f_1$.
Hence, solving \eqref{def L Strategy 2} is equivalent to solving 
\begin{equation*}
\mc{L}_0z=f_2,\
\end{equation*}
where 
\begin{gather}\label{def L02}
 \mc{L}_0z:=
 \partial_tz_2-\Div (d_2\nabla z_2)
 -\partial_{x_1} z_1
 - g_{22}\cdot\nabla z_2- a_{22}z_2\mbox{ in } \omega_T.
\end{gather}
Hence, finding a differential operator   $\mc{M}$ such that \eqref{strat:LMB} is satisfied 
is now equivalent to finding a differential operator   $\mc{M}_0$ such that 
\begin{equation}\label{L0M0}
 \mc{L}_0\circ\mc{M}_0=Id.
\end{equation}
%
%
%
We can remark that equality \eqref{L0M0} is formally equivalent to
\begin{equation}\label{M*L*}
\mc{M}_0^*\circ  \mc{L}_0^*=Id,
\end{equation}
where the formal adjoint $ \mc{L}_0^*$ of the operator $ \mc{L}_0$ is given 
for every $\varphi\in \mc{C}^{\infty}(\overline{\omega}_T)$ by
\begin{equation}\label{def L*}
\begin{array}{rcl}
 \mc{L}^*_0\varphi&:=& \left(\begin{array}{c}
 {\mc{L}}_1\varphi\\\noalign{\smallskip}
 {\mc{L}}_{2}\varphi
        \end{array}\right)
= \left(\begin{array}{c}\partial_{x_1}\varphi \\\noalign{\smallskip}
-\partial_t(\varphi)-\Div ( d_2\nabla (\varphi))
+\Div (g_{22} \varphi)
- a_{22}\varphi\\\noalign{\smallskip}
        \end{array}\right).
\end{array}\end{equation}
Operator $\mc{L}_2$ can be rewritten as 
\begin{equation*}
\begin{array}{rcl}
 {\mc{L}}_{2}\varphi
=-\partial_t\varphi-\sum\limits_{i,j=1}^Nd_2^{ij}\partial_{x_ix_j}\varphi
+ \sum\limits_{i=1}^N\widetilde{g}_{22}^i\partial_{x_i}\varphi
+\widetilde{a}_{22}\varphi,
\end{array}\end{equation*}
where $\widetilde{g}_{22}^i$ and $\widetilde{a}_{22}$ 
are given in  \eqref{resol alg:def g22 tilde}. 
Let us first consider the following linear combination of $\mc{L}_1$ and  $\mc{L}_2$:
\begin{equation*}
\begin{array}{rcl}
 {\mc{L}}_{3}\varphi
 &=& {\mc{L}}_{2}\varphi
 -[-2\sum\limits_{i=2}^Nd_2^{i1}\partial_{x_i}
 + \widetilde{g}_{22}^1]\mc{L}_1\varphi\\
&=&-\partial_t\varphi-\sum\limits_{i,j=2}^Nd_2^{ij}\partial_{x_ix_j}\varphi
+ \sum\limits_{i=2}^N\widetilde{g}_{22}^i\partial_{x_i}\varphi
+\widetilde{a}_{22}\varphi,
\end{array}\end{equation*}
Lemma \ref{lemme:resol alg} leads to the algebraic resolvability of 
System \eqref{def L Strategy 2} under Condition \ref{cond:modul}.

Concerning the order of $\mathcal M$,
 if we follow the proof of Lemma \ref{lemme:resol alg} step by step,
we apply at most $N\times(N-1)/2$  operators of order two to eliminate the terms $d_2^{ij}\partial_{x_ix_j}$ with $i,j\in\{2,...,N\}$ (thanks to the symmetry property of $d_2$), 
then at most $N-1$  operators of order one
for the term $\widetilde{g}_{22}^i\partial_{x_i}$ with $i\in\{2,...,N\}$
and finally an operator of order at most one for $\partial_t$.
Thus the operator $\mathcal M$ is of order  at most $N\times(N-1)+(N-1)+1=N^2.$
\cqfd

We are now ready for the proof of Theorem \ref{theo:positif}.

\textbf{Proof of Theorem \ref{theo:positif}.}\\
We apply Proposition \ref{prop:contr regul} with $k=N^2+1$ and obtain the existence of  two constants $K>0$ and $C>0$ 
such that for every initial condition $y^0\in L^2(\Omega)^2$ 
one can find a control $\widehat u\in \mathcal{C}^{N^2+1}(\overline{Q_T})$  verifying 
$\Supp(\widehat{u})\subset \subset \omega_T$ for which 
the solution $\widehat y$  to System \eqref{strat:syst lin pb ana} is equal to zero at time $T$
and the following estimate holds:
\begin{equation}\label{regexp}
\|\widehat u\|_{\mathcal{C}^{N^2+1}(Q_T)^2}\leqslant C_k\|y^0\|_{L^2(\Omega)^2}.
\end{equation}

Now, using Proposition \ref{prop:resol alg}, locally on $\omega_T$ there exists a solution 
$(z,v)\in \mc{C}^1(\overline{Q}_T)^3\subset W(0,T)^2\times L^2(Q_T)$ to the following control problem:
\begin{equation*}
 \left\{\begin{array}{ll}
\partial_tz_1=\Div (d_1\nabla z_1)+g_{11}\cdot\nabla z_1+g_{12}\cdot\nabla z_2+a_{11}z_1+a_{12}z_2+\widehat{u}_1+v&\mbox{in }\omega_T,\\\noalign{\smallskip}
\partial_tz_2=\Div (d_2\nabla z_2)+\partial_{x_1} z_1+g_{22}\cdot\nabla z_2+a_{22}z_2+\widehat{u}_2&\mbox{in } \omega_T,
        \end{array}
\right.
\end{equation*}
with $(\widehat u_1,\widehat u_2):=\widehat u$. Moreover, since $\Supp(z)\subset \subset \omega_T$, we have 
$z(0,\cdot)=z(T,\cdot)=0$ in $\Omega$.

We conclude by remarking  that   $(y,u):=(\widehat{y}-z,-v)$ is a solution to 
System (\ref{system primmal}) which
satisfies $y(T,\cdot)\equiv0$ in $\Omega$.

\cqfd

\section{Proof of Theorem \ref{theo:negatif2}}

\hspace*{4mm} 

Let $\omega_1$  be a nonempty regular open set satisfying $\omega\subset\subset \omega_1\subset\subset\Omega$. Let
$\theta$ be  a function of $\mathcal{C}^{\infty}(\overline{\Omega})$ satisfying
\begin{equation*}\left\{\begin{array}{l}
\theta=1\mbox{ in }\omega_0,\\\noalign{\smallskip}
\Supp(\theta)\subset\\omega,\\\noalign{\smallskip}
\theta>0\mbox{ in } \omega_1
\end{array}\right.
\end{equation*}
Consider the following system 
\begin{equation}\label{syst simplp}
\left\{\begin{array}{ll}
\partial_ty_1=\Delta y+\mathds{1}_{\omega}u&\mbox{in } Q_T,\\\noalign{\smallskip}
\partial_ty_2=\Delta y_2+ay_2+\partial_{x_1}(\theta y_1)&\mbox{in } Q_T,\\\noalign{\smallskip}
y=0&\mbox{on } \partial\Omega,\\\noalign{\smallskip}
y(0,\cdot)=y^0&\mbox{in }\Omega,
\end{array}\right.
\end{equation}
where $u\in L^2(Q_T)$ is the control and $a\in L^{\infty}(\Omega)$ will be specified later.
If we can control approximately  System \eqref{syst simplp}, 
then it implies that we are also able to control approximately the following equation:
\begin{equation}\label{syst primal}
\left\{\begin{array}{ll}
\partial_tz=\Delta z+az+\partial_{x_1}(\theta v)&\mbox{in } Q_T,\\\noalign{\smallskip}
z=0&\mbox{on } \partial\Omega,\\\noalign{\smallskip}
z(0,\cdot)=y^0_2&\mbox{in }\Omega,
\end{array}\right.
\end{equation}
where $v\in L^2((0,T),H^1(\Omega))$ is the control.
Since $\theta>0$ on $\omega_1$, the approximate controllability on the time interval $(0,T)$  of System \eqref{syst simplp} is equivalent to the following property, called the Fattorini criterion (see \cite[Theorem 1 \& Section 3]{olive_bound_appr_2014}):
\begin{theo}\label{theo:fattorini2}
System \eqref{syst primal} is approximately controllable on the time interval $(0,T)$, if and only if for every
$s\in\mathbb{C}$ and every $\varphi \in D (\Delta )$, we have
\begin{equation*}
\left. \begin{array}{ll}
-\Delta  \varphi-a\varphi  = s\varphi &\mbox{ in } \Omega\\\noalign{\smallskip}
\partial_{x_1} \varphi= 0&\mbox{  in }\omega_1
\end{array}\right\}
\Rightarrow \varphi = 0.
\end{equation*}
\end{theo}
Since $\omega_1\subset\subset\Omega$,
Then there exists a open set $\omega_2$ such that $\omega_1\subset\subset \omega_2\subset\subset \Omega$.
The first eigenfunction $\varphi_1$ of $-\Delta$ is well-known to be positive in $\Omega$, 
so we can define a function $\varphi\in \mathcal{C}^{\infty}(\overline\Omega)$ satisfying
\begin{equation*}
\left\{\begin{array}{ll}
\varphi =\varphi_1 &\mbox{ in } \Omega\backslash \omega_2,\\\noalign{\smallskip}
\varphi = 1&\mbox{ in }\omega_1,\\\noalign{\smallskip}
\varphi>\delta>0&\mbox{ in }\omega_2.
\end{array}\right.
\end{equation*}
For instance, if $\Omega:=(0,\pi)$ and $\omega_1:=(2\pi/5,3\pi/5)$, as in Figure 1,
we may construct a function $\varphi\in \mathcal{C}^2([0,\pi])$  satisfying
\begin{equation*}
\left\{\begin{array}{ll}
\varphi(x) = \sin(x) &\mbox{ for every } x\in [0,\pi/5]\cup[4\pi/5,\pi ],\\\noalign{\smallskip}
\varphi(x) = 1&\mbox{ for every }x\in [2\pi/5,3\pi/5],\\\noalign{\smallskip}
\varphi>\delta>0&\mbox{ in }[\pi/5,4\pi/5].
\end{array}\right.
\end{equation*}
\begin{figure}[H]\begin{center}
\includegraphics[scale=0.9]{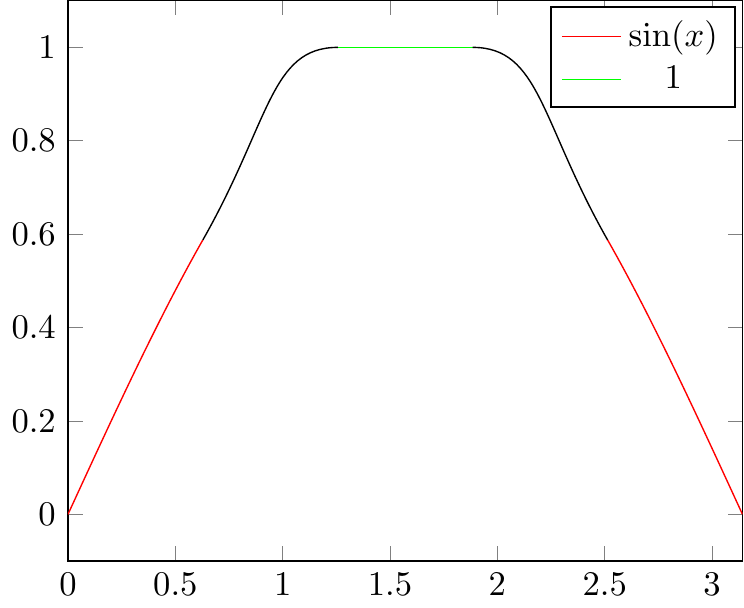}
\caption{Example of function $\varphi$ on $[0,\pi]$}
\end{center}\end{figure}
Consider
\begin{equation*}
a:=\dfrac{-\Delta \varphi-\varphi}{\varphi}.
\end{equation*}
Thanks to the definition of $\varphi$, is well defined in $\overline{\Omega}$ 
and is an element of $\mc{C}^{\infty}(\overline{\Omega})$. 
Thus $\varphi$ satisfies
\begin{equation*}
\left\{ \begin{array}{ll}
-\Delta  \varphi-a\varphi  = \varphi &\mbox{ in } \Omega,\\\noalign{\smallskip}
\partial_{x_1} \varphi = 0&\mbox{  in }\omega,\\\noalign{\smallskip}
\varphi\neq0.
\end{array}\right.
\end{equation*}
Using Theorem \ref{theo:fattorini2}, System \eqref{syst primal} is not approximately controllable on the time interval $(0,T)$.

\cqfd
\begin{rem}
Let us emphasize that in this case, as expected, Condition \ref{cond:modul} is not verified: on $\omega$ we have by definition $a_{22}=-1$, $g_{22}=0$ and $d^{ii}_2=0$ for every $i\in\{2,\ldots,N\}$, which implies that  $\tilde a_{22}=-1$ on $\omega$ and $\widetilde g_{22}=0$, hence
\begin{equation*}
\left\{\begin{array}{l}
\widetilde{a}_{22}\mbox{ is an element of the }
\mc{C}^0_{t,x_2,...,x_N}(\overline{\omega}_T)\mbox{-module }\\\noalign{\smallskip}
\left\langle 1,\widetilde{g}_{22}^2,...,\widetilde{g}_{22}^N,d_2^{22},...,d_2^{NN}\right\rangle_{\mc{C}^0_{t,x_2,...,x_N}(\overline{\omega}_T)}.
\end{array}\right.
\end{equation*}
This will also be the case for the potential constructed in the first part of the proof of Theorem \ref{theo:negatif}.
\end{rem}

\section{Proof of Theorem \ref{theo:negatif}}

\hspace*{4mm} 

Let $\Omega:=(0,\pi)$ and $\omega:=(7\pi/15,8\pi/15)$. 
Consider the following system 
\begin{equation}\label{syst simpl2}
\left\{\begin{array}{ll}
\partial_ty_1=\Delta y+\mathds{1}_{\omega}u&\mbox{in } Q_T,\\\noalign{\smallskip}
\partial_ty_2=\Delta y_2+ay_2+\partial_{x} y_1&\mbox{in } Q_T,\\\noalign{\smallskip}
y=0&\mbox{on } \Sigma_T,\\\noalign{\smallskip}
y(0,\cdot)=y^0&\mbox{in }\Omega,
\end{array}\right.
\end{equation}
where $u\in L^2(Q_T)$ is the control and $a\in C^{\infty}(\overline\Omega)$ will be specified later.

As in the previous section, it is well-known that the approximate controllability on the time interval $(0,T)$  of System \eqref{syst primal} is equivalent to the following property:

\begin{theo}\label{theo:fattorini}
System \eqref{syst simpl2} is approximately controllable on the time interval $(0,T)$, if and only if for every
$s\in\mathbb{C}$ and every $\varphi \in D (\Delta )$, we have
\begin{equation*}
\left. \begin{array}{ll}
-\Delta  \varphi-\partial_{x} \psi  = s\varphi &\mbox{ in } \Omega\\\noalign{\smallskip}
-\Delta\psi -a\psi= s\psi&\mbox{ in } \Omega\\\noalign{\smallskip}
\varphi= 0&\mbox{  in }\omega
\end{array}\right\}
\Rightarrow (\varphi,\psi) = (0,0).
\end{equation*}
\end{theo}

\begin{figure}[H]\begin{center}
\includegraphics[scale=0.9]{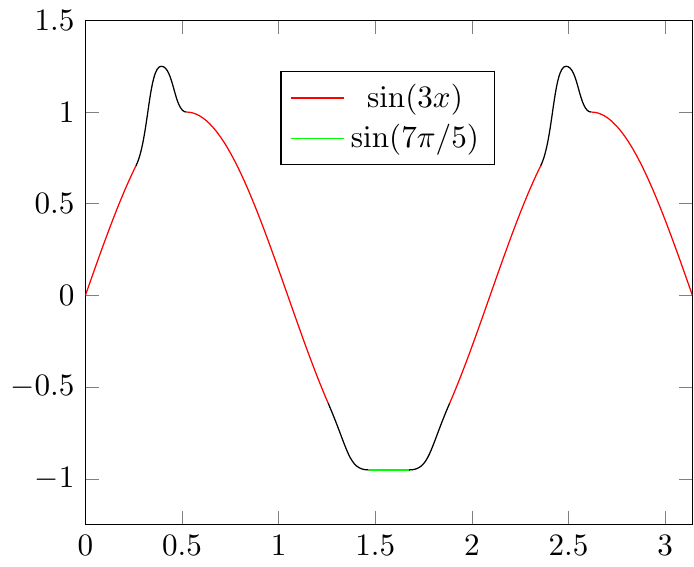}
\includegraphics[scale=0.9]{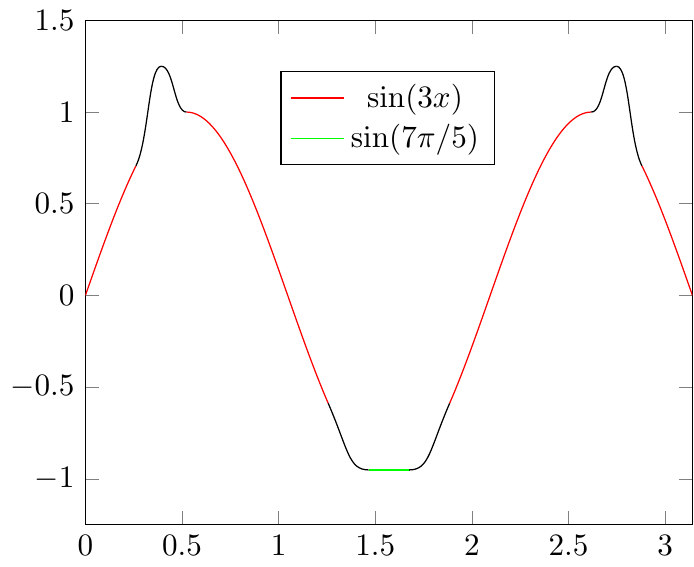}
\caption{Examples of function $\psi$ on $[0,\pi]$}\label{fig:contre exemple}
\end{center}\end{figure}

Let us construct three functions $\varphi$, $\psi$, $a\in\mc{C}^{\infty}(\overline{\Omega})$ satisfying
\begin{equation}\label{cond:constr simpl}
\left\{ \begin{array}{ll}
-\Delta  \varphi-\partial_{x} \psi  = 9\varphi &\mbox{ in } \Omega,\\\noalign{\smallskip}
-\Delta\psi -a\psi= 9\psi&\mbox{ in } \Omega,\\\noalign{\smallskip}
\varphi(0)=\varphi(\pi)=\psi(0)=\psi(\pi)=0,&\\\noalign{\smallskip}
\varphi= 0&\mbox{  in }\omega,\\\noalign{\smallskip}
\varphi\neq 0,~\psi\neq0&\mbox{  in }\Omega.
\end{array}\right.
\end{equation}
The idea will be to construct the function $\psi$  as a perturbation of 
$x\mapsto\sin(3x)$.
Consider $\psi$ a function of $\mc{C}^{\infty}(\overline{\Omega})\cap D(\Delta)$ satisfying
\begin{equation}\label{cond:constr simpl2}
\left\{
\begin{array}{ll}
\psi(x)=\sin(3x)+C_1\theta_1+C_2\theta_2+C_3\theta_3&\mbox{ for all }x\in\overline{\Omega},\\\noalign{\smallskip}
\psi(x)=\sin(7\pi/5)&\mbox{ for all }x\in\overline{\omega},\\\noalign{\smallskip}
|\psi(x) - \sin(3x)|<\varepsilon &\mbox{ for all }x\in[6\pi/15,7\pi/15]\cup[8\pi/15,9\pi/15],
\end{array}\right.
\end{equation}
where  $\theta_1,~\theta_2,~\theta_3$ 
are three nontrivial functions of $\mc{C}^{\infty}(\overline{\Omega})$ satisfying
\begin{equation}
\left\{ \begin{array}{l}
\Supp(\theta_1)\subset(\pi/12,\pi/6),\\\noalign{\smallskip}
\Supp(\theta_2)\subset(9\pi/12,5\pi/6),\\\noalign{\smallskip}
\Supp(\theta_3)\subset(5\pi/6,11\pi/12),\\\noalign{\smallskip}
\theta_1,~\theta_2,~\theta_3\geqslant0\mbox{ in } \Omega,
\end{array}\right.
\end{equation}
$\varepsilon>0$ small enough and
$C_1,~C_2,~C_2$ are three positive constants to determined
(See Figure \ref{fig:contre exemple} for some examples of function $\psi$).
Let us remark that, for a constant $\alpha\in\mb{R}$ to determined,
the function $\varphi\in\mc{C}^{\infty}(\overline{\Omega})$ defined for all $x\in \overline{\Omega}$ by
\begin{equation*}
\begin{array}{rcl}
\varphi(x)&:=&\alpha\sin(3x)-\frac13\displaystyle\int_0^x\sin(3(x-y))\partial_x\psi(y)dy
\end{array}
\end{equation*}
is solution to the first equation of \eqref{cond:constr simpl}.
In order to apply Theorem \ref{theo:fattorini}, let us first prove that $C_1$ and $\alpha$ can be chosen such that $\varphi=0$ in $\omega$.
Since $\psi=\sin(7\pi/5)$ in $\omega$, 
\begin{equation*}
\begin{array}{rcl}
\varphi(x)
&=&\left[\alpha-\frac13\cos(7\pi/5)\sin(7\pi/5)-\displaystyle\int_0^{7\pi/15}\sin(3y)\psi(y)dy\right]\sin(3x)\\
&&\hspace*{2cm}
+\left[\frac13\sin(7\pi/5)^2-\displaystyle\int_0^{7\pi/15}\cos(3y)\psi(y)dy\right]\cos(3x),
\end{array}
\end{equation*}
for all $x\in\omega$.
Since $\cos(3x)>0,~\sin(3x)>0$ for all $x$ in  $(\pi/12,\pi/6)$ and 
\begin{equation*}
\frac13\sin(7\pi/5)^2-\displaystyle\int_0^{7\pi/15}\cos(3y)\sin(3y)dy>0,
\end{equation*}
then, according to the last line of \eqref{cond:constr simpl2}, for $\epsilon$ small enough,
it is possible to choose $C_1>0$ in order to obtain 
\begin{equation*}
\frac13\sin(7\pi/5)^2-\displaystyle\int_0^{7\pi/15}\cos(3y)\psi(y)dy=0.
\end{equation*}
Thus, for $\alpha$ given by
\begin{equation*}
\alpha:=\frac13\cos(7\pi/5)\sin(7\pi/5)+\displaystyle\int_0^{7\pi/15}\sin(3y)\psi(y)dy,
\end{equation*}
we obtain $\varphi=0$ in $\omega$.
By definition of $\varphi$, we have $\varphi(0)=0$. 
Let us now prove that for some appropriate $C_2$ and $C_3$, we have $\varphi(\pi)=0$.
We remark that
\begin{equation*}
\varphi(\pi)=\frac13\displaystyle\int_0^{\pi}\cos(3y)\psi(y)dy.
\end{equation*}
Let us distinguish two cases: 
\begin{enumerate}
\item If 
\begin{equation}\label{quantity phi(pi)}
\frac13\displaystyle\int_0^{2\pi/3}\cos(3y)\psi(y)dy+\frac13\displaystyle\int_{2\pi/3}^{\pi}\cos(3y)\sin(3y)dy
\end{equation}
is negative, then, using the fact that $\sin(3x),~\cos(3x)>0$ 
for all $x\in (9\pi/12,5\pi/6)$, one can choose  $C_3:=0$ and find some some $C_2>0$ such that
$\varphi(\pi)=0$. 
\item If now the quantity \eqref{quantity phi(pi)} is positive, 
since $\sin(3x)>0$ and $\cos(3x)<0$ for all $x\in (5\pi/6,11\pi/12)$, 
 one can choose  $C_2:=0$ and find some some $C_3>0$ such that
$\varphi(\pi)=0$. 
\end{enumerate}
The function $\psi$ will have one of the two following forms

To satisfy the second equality in \eqref{cond:constr simpl}, 
we define the function $a\in\mc{C}^{\infty}(\overline{\Omega})$ as follows
\begin{equation}\label{def:a}
a:=\dfrac{-\Delta \psi-9\psi}{\psi}.
\end{equation}
This function $a$ is bounded since 
at each point where $\psi$ is null, i.e. at $0$, $\pi/3$, $2\pi/3$ and $\pi$, 
there exists a neighbourhood  in which $\psi(x)$ is equal to $\sin(3x)$.
Thus the constructed $\varphi$, $\psi$ and $a$ verify  \eqref{cond:constr simpl}.
Using Theorem \ref{theo:fattorini}, System \eqref{syst simpl2} is not approximately controllable on the time interval $(0,T)$.

Let us now prove the second item of Theorem \ref{theo:negatif}.
We remark that it is possible to chose $\theta_1=\mbox{exp}$ in $\omega_1\subset(\pi/12,\pi/6)$ with $\omega_1$ small enough. 
Then $a$ is defined in $\omega_1$ for all $x\in\omega_1$ by
$$a(x)=\dfrac{-10C_1\exp(x)}{\sin(3x)+C_1\exp(x)}.$$
Thus $a$ satisfies Condition  \ref{cond:modul} for $\omega:=\omega_1$, that is $a$ is non-constant in the space variable on $\omega_1$. 
We conclude applying Theorem \ref{theo:positif} for $\omega:=\omega_1$.

\cqfd


%

\vspace{0.3cm}

\textbf{{\Large Funding}}

Pierre Lissy is partially supported by the project IFSMACS funded by the french Agence Nationale de la Recherche, 2015-2019 (Reference:  ANR-15-CE40-0010).

\vspace{0.3cm}

\textbf{{\Large Conflict of Interest}}

The authors declare that they have no conflict of interest.


\end{document}